\documentclass[envcountsame]{llncs}
\usepackage{amssymb}
\usepackage{amsmath}
\usepackage[margin=1.4in]{geometry}

\newcommand{\A}{\mathfrak{A}}
\newcommand{\B}{\mathcal{B}}

\newcommand{\F}{\mathcal{F}}
\newcommand{\M}{\mathcal{M}}
\newcommand{\N}{\mathbb{N}}
\newcommand{\Q}{\mathbb{Q}}
\newcommand{\U}{\mathcal{U}}

\newcommand{\X}{\mathcal{X}}

\newcommand{\dfc}{\mathbf{d}}
\newcommand{\ld}{\mathbf{ld}}
\newcommand{\score}{\mathbf{score}}
\newcommand{\WRONG}{\ensuremath{\textsc{Wrong}}}
\newcommand{\SUCC}{\ensuremath{\textsc{Succ}}} 
\newcommand{\BEST}{\ensuremath{\textsc{Best}}}

\newcommand{\Leb}{\lambda}

\renewcommand{\epsilon}{\varepsilon}

\newcommand{\K}{\mathrm{K}}
\newcommand{\fs}{2^{<\omega}}
\newcommand{\dist}{\rho}
\newcommand{\prefix}{\preceq}
\newcommand{\Mcs}{\mathcal{P}}

\newcommand{\cs}{2^\omega}
\newcommand{\uh}{{\upharpoonright}}
\newcommand{\leqx}{\preceq}
\newcommand{\emptystr}{\Lambda}

\newcommand{\md}{\mskip 1mu|\mskip 1mu}

{\itshape}{\rmfamily}
{\itshape}{\rmfamily}
{\itshape}{\rmfamily}
{\itshape}{\rmfamily}

\spnewtheorem{remarkbold}[theorem]{Remark}{\bfseries}{\itshape}

{\bfseries}{\rmfamily}
{\bfseries}{}

\begin{document}

\title{Algorithmic identification of probabilities is hard}
\author{Laurent Bienvenu\inst{1}
\and Santiago Figueira \inst{2}
\and Beno\^it Monin\inst{3}
\and Alexander Shen \inst{1}}

\institute{LIRMM, CNRS \& Universit\'e de Montpellier, France\\
\and
Universidad de Buenos Aires and CONICET, Argentina\\
\and 
LACL, Universit\'e Paris 12, France
}
\maketitle

\begin{abstract}
Suppose that we are given an infinite binary sequence which is random for a Bernoulli measure of parameter~$p$. By the law of large numbers, the frequency of zeros in the sequence tends to~$p$, and thus we can get better and better approximations of~$p$ as we read the sequence. We study in this paper a similar question, but from the viewpoint of inductive inference. We suppose now that $p$ is a computable real, and one asks for more: as we are reading more and more bits of our random sequence, we have to eventually guess the exact parameter~$p$ (in the form of its Turing code). Can one do such a thing uniformly for all sequences that are random for computable Bernoulli measures, or even for a `large enough' fraction of them? In this paper, we give a negative answer to this question. In fact, we prove a very general negative result which extends far beyond the class of Bernoulli measures. We do however provide a weak positive result, by showing that looking at a sequence~$X$ generated according to some computable probability measure, we can eventually guess a sequence of measures with respect to which~$X$ is random in Martin-L\"of's sense. 

\end{abstract}

\section{Introduction}

\subsection{Inductive inference}

The study of learnability of computable sequences is concerned with the following problem. Suppose we have a black box that generates some infinite computable sequence of bits $X=X(0) X(1) X(2),\ldots$ We do not know the program running in the box, and want to guess it by looking at finite prefixes
\[
X \uh n=X(0) \ldots X(n-1)
\]
for increasing values of~$n$. There could be different programs that produce the same sequence, and it is enough to guess one of them (since there is no way to distinguish between them by just looking at the output bits). The more bits we see, the more information we have about the sequence. For example, it is hard to say something about a sequence seeing only that its first bit is a~$1$, but looking at the prefix
\[
110010010000111111011010101000
\]
one may observe that this is a prefix of  the binary expansion of~$\pi$, and guess that the machine inside the box does exactly that (though the machine may as well produce the binary expansion of, say, $47627751/15160384$).

The hope is that, as we gain access to more and more bits, we will \emph{eventually} figure out how the sequence~$X$ is generated. More precisely, we hope to have a total computable function~$\A$ from strings to integers such that for every~computable $X$, the sequence
\[
\A(X \uh 1),\, \A(X \uh 2),\,  \A(X \uh 3),\ldots
\]
converges to a program (= index of a computable function) that computes $X$. This is referred to as~\emph{identification in the limit}, and can be understood in (at least) two ways. Indeed, assuming that we have a fixed effective enumeration $(\varphi_e)_{e \in \N}$ of partial computable functions from $\N$ to $\{0,1\}$, we can define two kinds of success for an algorithm~$\A$ on a computable sequence~$X$:

\begin{itemize}
\item Strong success: the sequence $e_n=\A(X \uh n)$ converges to a single value~$e$ such that $\varphi_e=X$ (i.e., $\varphi_e(k)=X(k)$ for all~$k$). 
\item Weak success: the sequence $e_n=\A(X \uh n)$ does not necessarily converge, but $\varphi_{e_n}=X$ for all sufficiently large~$n$.
\end{itemize}
Here we assume that $\A(X\uh n)$ is defined for all $n$ or at least for all sufficiently large $n$.

The strong type of success is often referred to as~\emph{explanatory} (EX), see, e.g., Definition~VII.5.25 in~\cite[p.~116]{OdifreddiVolume2}. The second type is referred (see Definition~VII.5.44, p.~131 in the same book) as~\emph{behaviorally correct} (BC). Note that it is obvious from the definition that strong success implies weak success.  

It would be nice to have an algorithm that succeeds on all computable sequences. However, it is impossible even for weak success: for every (total) algorithm~$\A$, there is a computable~$X$ such that $\A$ does not weakly succeed on~$X$. 
The main obstacle is that certain machines are not total (produce only finitely many bits), and distinguishing total machines from non-total ones cannot be done computably. 

However, some \emph{classes} of computable sequences can be learned, i.e., there exists a total algorithm that succeeds on all elements of the class. Consider for example the class of primitive recursive functions. This class can be effectively enumerated, i.e., there is a total computable function $f$ such that $(\varphi_{f(e)})_{e \in \N}$ is exactly the family of primitive recursive functions. Now consider the algorithm~$\A$ such that $\A(\sigma)$ returns the smallest $e$ such that $\varphi_{f(e)}(i)=\sigma(i)$ for all~$i < |\sigma|$ (such an $e$ always exists, since every string is a prefix of a primitive recursive sequence). It is easy to see that if $X$ is primitive recursive, $\A$ succeeds on $X$, even in the strong sense (EX).

The theory of learnability of computable sequences (or functions) is precisely about determining which classes of functions can be learned. This depends on the learning model, the type of success, of which there are many variants. We refer to the survey by Zeugman and Zilles~\cite{ZeugmannZ2008} and to~\cite[Chapter VII]{OdifreddiVolume2} for a panorama of the field. 

\subsection{Learning measures}

	Recently, Vit\'anyi and Chater~\cite{ChaterV2017} proposed to study a related problem. Suppose that instead of a sequence that has been produced by a deterministic machine, we are given a sequence that has been generated by a randomized algorithmic process, i.e., by a Turing machine that has access to a fair coin and produces some output sequence on the one-directional write-only output tape. The output sequence is therefore a random variable defined on the probabilistic space of fair coin tossings. We assume that this machine is \emph{almost total}.\footnote{This requirement may look unnecessary. Still the notion of algorithmic randomness needed for our formalization is well-defined only for computable measures, and machines that are not almost total may not define a computable measure.} This means that the generated sequence is infinite with probability~$1$.
	
Looking at the prefix of the sequence, we would like to guess which machine is producing it. For example,  for the sequence
\[
000111111110000110000000001111111111111
\]
we may guess that it has been generated via the following process: start with~$0$ and then choose each output bit to be equal to the previous one with probability, say, $4/5$ (so the change happens with probability $1/5$), making all the choices independently.\footnote{The probability $4/5$ is not a dyadic rational number, but still can be simulated by an almost total machine using a fair coin.}

So what should count as a good guess for some observed sequence? Again, there is no hope to distinguish between two  processes that have the same output distribution. So our goal should be to reconstruct the output distribution and not the specific machine.

But even this is too much to ask for. Assume that we have agreed that some machine $M$ with output distribution~$\mu$ is a plausible explanation for some sequence $X$. Consider another machine $M'$ that starts by tossing a coin and then (depending on the outcome) either generates an infinite sequence of zeros or simulates $M$. If $X$ is a plausible output of $M$, then $X$ is also a plausible output of $M'$, because it may happen (with probability $1/2$) that $M'$ simulates $M$.

A reasonable formalization of a `good guess' is provided by the theory of algorithmic randomness. As Chater and Vit\'anyi recall, there is a widely accepted formalization of ``plausible outputs'' for an almost total probabilistic machine with output distribution~$\mu$: the notion of Martin-L\"of random sequences with respect to~$\mu$. These are the sequences that pass all effective statistical tests for the measure~$\mu$, also known as $\mu$-Martin-L\"of tests. (We assume that the reader is familiar with algorithmic randomness and Kolmogorov complexity. The most useful references for our purposes are~\cite{Gacs2005} and~\cite{LiV2008}). Having this notion in mind, the natural way to extend learning theory to the probabilistic case is as follows: 

\begin{quote}
\emph{A class of computable measures~$\M$ is learnable if there exists a total algorithm $\A$ such that for every sequence~$X$ that is Martin-L\"of random for some measure in $\M$, the sequence $$\A(X \uh 1), \A(X \uh 2), \A(X \uh 3),\ldots$$ identifies in the limit a measure $\mu \in \M$ such that $X$ is Martin-L\"of random with respect to $\mu$.}
\end{quote}

Like in the classical case, there are several ways one can interpret the notion of `identifying in the limit. We will come back to this after having introduced some basic notation and terminology related to computable measures (for now one may think of a computable measure as an output distribution of an almost total probabilistic machine).

\subsection{Background and notation}

We denote by $\cs$ the set of infinite binary sequences and by $\fs$ the set of finite binary sequences (or \emph{strings}). The length of a string $\sigma$ is denoted by~$|\sigma|$. The empty string (string of length~$0$) is denoted by $\emptystr$. For two strings $\sigma$, $\tau$ we write $\sigma\prefix \tau$ if $\sigma$ is a prefix of $\tau$. The $n$-th element of a sequence $X(0)X(1)\ldots$ is the value $X(n-1)$ (assuming that the length of $X$ is at least $n$); the string $X\uh n = X(0)X(1) \ldots X(n-1)$ is the \emph{$n$-bit prefix of $X$}.  We write $\sigma \prefix X$ if the string $\sigma$ is a prefix of the infinite sequence~$X$ (i.e., $X \uh |\sigma| = \sigma$).
The space $\cs$ is endowed with the distance~$d$ defined by
\[
d(X,Y) = 2^{- \min \{n: X(n) \ne Y(n)\}}.
\]
This distance is compatible with the product topology generated by \emph{cylinders}
\[
[\sigma] = \{X \in \cs \; : \; \sigma \prefix X\}.
\]
A cylinder is both open and closed (= \emph{clopen}). Thus, any finite union of cylinders is also clopen. It is easy to see, by compactness, that the converse holds: every clopen subset of~$\cs$ is a finite union of cylinders. We say that a clopen set~$C$ has \emph{granularity at most~$n$} if $C$ is a finite union of some cylinders $[\sigma]$ with $|\sigma|=n$. We denote by $\Gamma_n$ the family of clopen sets of granularity at most~$n$.

We now give a brief review of the ``computable analysis" aspects of the space of probability measures. For a more thorough exposition of the subject, the main reference is~\cite{Gacs2005}. 

The space of Borel probability measures over $\cs$ is denoted by $\Mcs$. In the rest of the paper, when we talk about a `measure', we mean an element of the space~$\Mcs$. This space is equipped with the weak topology, that is, the weakest topology such that for every $\sigma$, the application $\mu \mapsto \mu([\sigma])$ is continuous as a function from $\Mcs$ to $\mathbb{R}$. Several classical distances are compatible with this topology; for example, one may use the distance~$\dist$ constructed as follows. For $\mu,\nu \in \Mcs$, let $\dist_n(\mu,\nu)$ (for an integer~$n$) be the quantity
\[
\dist_n(\mu,\nu)=\max_{C \in \Gamma_n} |\mu(C)-\nu(C)|
\]
and then set
\[
\dist(\mu,\nu) = \sum_n 2^{-n} \dist_n (\mu,\nu).
\]
The \emph{open} (resp.\ \emph{closed}) \emph{ball $\B$ of center $\mu$ and radius~$r$} is the set of measures~$\nu$ such that $\dist(\mu,\nu)<r$ (resp.\ $\dist(\mu,\nu) \leq r$). In the space of measures, the closure $\overline{\B}$ of the open ball $\B$ of center $\mu$ and radius~$r$ is the closed ball of center $\mu$ and radius~$r$.  

The space $\Mcs$ is separable, i.e., has a countable dense set of points. An easily describable one is the set $\mathcal{I}$ consisting of measures $\{\delta_\sigma\}_{\sigma \in \fs}$, where $\delta_\sigma$ is the Dirac measure concentrated on the point $\sigma 0^\omega$, and all rational convex combinations of such measures. Note that every member of $\mathcal{I}$ has a finite description: it suffices to give the list of $\sigma$'s together with the rational coefficients of the linear combination. Thus one can safely talk about computable functions from/to $\mathcal{I}$.

The set $\mathcal{I}$, together with the distance~$\rho$, make $\Mcs$ a computable metric space~\cite{Gacs2005}. Each point~$\mu \in \Mcs$ can be written as the limit of a sequence $(q_1, q_2, \ldots)$ of points in~$\mathcal{I}$ where $\rho(q_i,q_j) \leq 2^{-i}$ for $i<j$. Such a sequence is called a \emph{fast Cauchy name} for~$\mu$. We say that a measure $\mu$ is \emph{computable} if there is a total computable function $\varphi_e: \N \rightarrow \mathcal{I}$ such that $(\varphi_e(n))_{n \in \N}$ is a fast Cauchy name for~�$\mu$. Such an~$e$ is called an \emph{index for $\mu$}. 

At this point the way we view measures ---  as points of the space $\Mcs$ --- does not match the presentation of the introduction, where we asked the learning algorithm to guess, on prefixes of input~$X$, a sequence of probabilistic machines~$M_i$ such that for almost all~$i$, the machine $M_i$ is almost total and $X$ is a plausible output for $M_i$. The reason is that in fact there are three ways one can think of measures, which are equivalent for our purposes:

\begin{itemize}
\item[(a)] A measure is a point of $\Mcs$. 
\item[(b)] By Caratheodory's theorem, a measure $\mu$ can be identified with the function $\sigma \mapsto \mu([\sigma])$: for every function $f: \fs \rightarrow [0,1]$ such that $f(\emptystr)=1$ and $f(\sigma 0)+f(\sigma 1)=f(\sigma)$ there is a unique measure $\mu$ such that $\mu([\sigma])=f(\sigma)$ for all $\sigma$. For example, the uniform measure $\Leb$ is the unique measure such that $\Leb([\sigma])=2^{-|\sigma|}$ for all $\sigma$, and the Bernoulli measure $\beta_p$ of parameter $p \in [0,1]$ is the unique measure satisfying $\beta_p ([\sigma 1]) = p \cdot \beta_p([\sigma])$ for all~$\sigma$. 

\item[(c)] Consider a Turing functional~$M$, which one might think of as a Turing machine with a read-only input tape, a work tape and a write-only output tape. We say that $M$ is defined on~$X$ if $M$ prints an infinite sequence $Y$ on the output tape given~$X$ on the input tape. When $M$ is defined on $\Leb$-almost every~$X$, where $\Leb$ is the uniform Lebesgue measure on a Cantor space that corresponds to the fair coin tossings, we say that $M$ is \emph{almost total}. Then the function 
\[
\mu_M(\sigma) = \Leb \{X \colon M(X) \succeq \sigma\} 
\]
defines a measure in the sense of item (b). This measure corresponds to the distribution of a random variable that is the output of $M$ on the sequence of uniform independent random bits.
 
\end{itemize}

These approaches are equivalent both in the classical and effective realm, as is well known. The corresponding classes of measures coincide; moreover, one can computably convert an algorithm representing a computable measure according to one of the definition, into other representations. However, depending on the context one characterization may be much easier to handle than the others. And indeed, the techniques of the next section where we will prove our main negative result are of analytic nature, so characterization (a) will be more convenient, while the positive result of the last section has a more `algorithmic flavor', for which characterization~(c) will be better suited. 

The \emph{randomness deficiency}\footnote{This version of randomness deficiency function is sometimes called ``uniform probability-bounded randomness deficiency''; however, we do not use the other versions and call it just ``randomness deficiency''.} function $\dfc$ is the largest, up to additive constant, function $f : \cs \times \Mcs \rightarrow \N \cup \{\infty\}$ such that
\begin{itemize}
\item $f$ is lower semi-computable (i.e., $f^{-1}((k,\infty])$ is an effectively open\footnote{An effectively open set is a union of a computably enumerable set of rational balls (or products of balls, since we consider a product space).} subset of the product space $\cs \times \Mcs$, uniformly in~$k$);
\item for every $\mu \in \Mcs$, for every integer~$k$, the inequality $\mu \{X\colon f(X,\mu) > k\} < 2^{-k}$ holds.
\end{itemize}

We use the usual notation $\dfc(X \md \mu)$ instead of $\dfc(X,\mu)$. We say that $X$ is (uniformly) random relative to measure $\mu$ if $\dfc(X \md \mu) < \infty$. For computable measures this notion coincides with the classical notion of Martin-L\"of randomness.

We end this introduction with a discussion on a concept we will need to state the main theorem of Section~\ref{sec:main-result}: orthogonality. Two measures $\mu, \nu \in \Mcs$ are said to be \emph{orthogonal} if there is a Borel set $\X \subseteq \cs$ such that $\mu(\X)=1$ and $\nu(\X)=0$ (taking the complement of~$\X$, we see that  ortogonality is a symmetric relation). This is equivalent to the following condition: for each $\epsilon>0$ there is a set $\X_\epsilon$ such that $\mu(\X_\epsilon)\geq 1-\epsilon$ and $\nu(\X_\epsilon) < \epsilon$ (indeed, one can then take $\X = \bigcap_i \bigcup_j \X_{2^{-i-j}}$). 

The class of Bernoulli measures provides an easy example of orthogonality: if $p \ne q$, the Bernoulli measures $\beta_p$ and $\beta_q$ (see  the definition above) are orthogonal (by the law of large numbers, taking for $\X$ the set of sequences with a limit frequency of ones equal to $p$, we have $\beta_p(\X)=1$ and $\beta_q(\X)=0$). 

The important fact we need is that when two computable measures $\mu$ and $\nu$ are orthogonal, they share no random element, i.e, $\dfc(X\md \mu)$ and $\dfc(X\md \nu)$ cannot both be finite for any $X$. For a proof of this result, see for example~\cite{BienvenuM2009}. 

\subsection{Learning models}

Most classical learning models for computable sequences can be adapted to our probabilistic setting. For example, the EX and BC models mentioned have the following natural counterparts (we give them the same names, as this should create no confusion). 

\begin{definition}
Let $X \in \cs$ and $\A: \fs \rightarrow \N$ a total algorithm. We say that:
\begin{itemize}
\item $\A$ EX-succeeds on~$X$ if $\A(X \uh n)$ converges to a value~$e$ that is an index for a computable measure~$\mu$ with respect to which~$X$ is Martin-L\"of random. 
\item $\A$ BC-succeeds on~$X$ if there exists a computable measure~$\mu$ such that for almost all~$n$, $\A(X \uh n)$ is an index for~$\mu$ and $X$ is Martin-L\"of random with respect to~$\mu$. 
\end{itemize}
\end{definition}

There are also some natural learning models we can define that are more specific to the probabilistic setting. As we discussed, for a given~$X$ that is Martin-L\"of random with respect to some computable measure, there are several (actually, infinitely many) computable measures with respect to which~$X$ is Martin-L\"of random.  Thus we could allow the learner to propose different measures at each step and not converge to a specific measure, as long as almost all of them are good explanations for the observed~$X$. To measure how good an explanation is, we use the randomness deficiency, thus it makes sense to make the distinction between learning with bounded randomness deficiency and with unbounded randomness deficiency. 

\begin{definition}
Let $X \in \cs$ and let $\A: \fs \rightarrow \N$ be a total algorithm. We say that:
\begin{itemize}
\item $\A$ BD-succeeds on~$X$ if there exists a constant~$d$ such that for almost all~$n$, $\A(X \uh n)$ is an index for a computable measure with respect to which~$X$ is Martin-L\"of random, with randomness deficiency at most~$d$. 
\item $\A$ UD-succeeds on~$X$ if for almost all~$n$, $\A(X \uh n)$ is an index for a computable measure with respect to which~$X$ is Martin-L\"of random. 
\end{itemize}
\end{definition}
(`BD' and `UD' stand for `bounded deficiency' and `unbounded deficiency'). Our four learning models are by no means an exhaustive list of possibilities. Just like the classical learning theory offers a wide variety of models, one could define a wealth of alternative models (partial learning, team learning, etc.) in our setting. This would take us far beyond the scope of the present paper and we leave this for further investigation. 

Let us note in passing that the four learning models we have presented form a hierarchy, namely:
\[
\text{EX-success} \ \Rightarrow \ \text{BC-success} \ \Rightarrow \ \text{BD-success} \ \Rightarrow \ \text{UD-success}  
\]

The fact that EX-success implies BC-success and that BD-success implies UD-success is immediate from the definition. To see that BC-success implies BD-success, recall that in our definition the randomness deficiency depends only on the measure but not on the algorithm that computes it. So if the learning algorithm BC-succeeds on some sequence $X$, i.e., outputs the same measure (its code) for all sufficiently large prefixes of $X$, then the deficiency of $X$ with respect to this measure will be a constant and therefore the algorithm BD-succeds on $X$.

\section{Identifying measures is hard}\label{sec:main-result}

Now that we have given a precise definition of various learning models for computable probability measures, there are some obvious questions we need to address, the first of which is: For each of the above learning models, is there a single algorithm~$\A$ that  succeeds on all sequences~$X$ that are random with respect to some computable measure? This measure can be different for different $X$.) And if not, are there natural classes of measures for which there is an algorithm which succeeds on all~$X$ that are random with respect to some measure in this class? 

The starting point of this paper was a claim made in a preprint of Vit\'anyi and Chater~\cite{ChaterV-false}, where it was stated that there exists an algorithm~$\A$ that EX-succeeds on every~$X$ that is Martin-L\"of random with respect to a Bernoulli measure $\beta_p$ for some computable~$p$ (different for different $X$). Our results (Theorems~\ref{thm:main-weak} and~\ref{thm:main-strong})  imply that there is in fact no such algorithm. Vit\'anyi and Chater later corrected this claim and proved the following weaker statement. 

\begin{theorem}[Vit\'anyi--Chater~\cite{ChaterV2017}]
Let $(p_e)$ be a partial enumeration of computable reals in~$[0,1]$. If $E \subseteq \N$ is c.e.\ or co-c.e., and for all~$e \in E$, $p_e$ is defined, then there exists an algorithm~$\A$ that EX-succeeds on every~$X$ that is random with respect to some $\beta_{p_e}$ for some~$e \in E$.
\end{theorem}

This result implies, for example, that there is an algorithm that EX-succeeds on all~$X$ that are random with respect to some $\beta_q$ with $q$ a rational number.

We prove that this result cannot be extended to all computable parameters~$p$:

\begin{theorem}\label{thm:main-weak}
No algorithm~$\A$ can BD-succeed on every sequence~$X$ that is random with respect to some Bernoulli measure $\beta_p$ for computable~$p$. A fortiori, there is no algorithm~$\A$ can BD-succeeds on every sequence~$X$ that is random with respect to some computable measure.
\end{theorem}

We will in fact prove a more general theorem, replacing the class of Bernoulli measures by any class of measures having some ``reasonable" structural properties, and allowing the learning algorithm to succeed on a fraction of sequences only. 

\begin{theorem}\label{thm:main-strong}
Let $\M$ be  a subset of $\Mcs$ with the following properties\textup:
\begin{itemize}
\item $\M$ is effectively closed, i.e., its complement is effectively open: one can enumerate a sequence of rational open balls in $\Mcs$ whose union is the complement of $\M$.
\item $\M$ is computably enumerable, i.e., one can enumerate all rational open balls in $\Mcs$ that intersect~$\M$.
\item for every computable measure $\nu$, and every non-empty open subset of~$\M$ \textup(i.e., a non-empty intersection of an open set in $\Mcs$ with $\M$\textup) there is a computable $\mu$ in this open subset that is orthogonal to $\nu$.
\end{itemize}
Let also $\delta$ be a positive number. Then there is no algorithm~$\A$ such that for every computable $\mu \in \M$, the $\mu$-measure of sequences~$X$ on which $\A$ BD-succeeds is at least~$\delta$.
\end{theorem}

The notion of a computably (= recursively) enumerable closed set is standard in computable analysis, see~\cite[Definition 5.1.1]{Weihrauch2000}.

Note that the hypotheses on the class $\M$ are not very restrictive: many standard classes of probability measures have these properties. In particular, the class $\{\beta_p \colon p \in [0,1]\}$ of Bernoulli measures is such a class, which is why. So we get Theorem~\ref{thm:main-weak} as a corollary: there is no algorithm that can learn all Bernoulli measures (not to speak about all Markov chains). To see that the third condition is true for the class of Bernoulli measures, note that only countably many Bernoulli measures may be non-orthogonal to a given measure $\mu$: the sets $L_p$ of sequences with limit frequency $p$ are disjoint, so only countably many of them may have positive $\mu$-measure. It remains to note that every open non-empty subset of the class of Bernoulli measures has the cardinality of the continuum.

Let us give another example (beyond Bernoulli measures and Markov chains) that satisfies the requirements of Theorem~\ref{thm:main-strong}. In this example, the probability of the next bit to be $1$ may depend on many of the previous bits. For every parameter $p \in [0,1]$, consider the measure $\mu_p$ associated to the stochastic process that generates a binary sequence bit by bit as follows: the first bit is $1$, and the conditional probability of $1$ after $\sigma 1 0^k$ is $p/(k+1)$. One can check that the class $\Mcs=\{\mu_p: p \in [0,1]\}$ satisfies the hypotheses of the theorem (observe that~$p$ can easily be reconstructed from the sequence that is random with respect to $\mu_p$).

Note also that these hypotheses are not added just for convenience: although they might not be optimal, they cannot be outright removed. If we do not require the class $\M$ to be effectively closed, compactness, then the class of Bernoulli measures $\beta_p$ with \emph{rational} parameter $p$ would qualify, but Vit\'anyi and Chater's theorem tells us that there \emph{is} an algorithm that correctly identifies each of the measures in the class with probability~$1$. The third condition is important, too. Consider the measures $\beta_0$ and $\beta_1$ concentrated on the sequences~$0000\ldots$ and $1111\ldots$ respectively. Then the class $\M= \{p \beta_0 + (1-p)\beta_1  \; : \; p \in [0,1]\}$ is indeed effectively closed and computably enumerable, but it is obvious that there is an algorithm that succeeds with probability~$1$ for all measures of that class (in the strongest sense: the first bit determines the entire sequence). For the second condition we do not have a counterexample showing that it is really needed, but it is true for all the natural classes (and it is guaranteed to be true if $\M$ has a computable dense sequence). 

\medskip
The rest of this section is devoted to the proof of Theorem~\ref{thm:main-strong}. 
\medskip

Fix a subset $\M$ of $\Mcs$ satisfying the hypotheses of the theorem, and some $\delta >0$. Assume for the sake of contradiction that there is a total algorithm $\A$ such that for every computable $\mu \in \M$, the $\mu$-measure of sequences~$X$ on which $\A$ BD-succeeds is at least~$\delta$. In the rest of the proof, by ``success'' we always mean BD-success.

We may assume without loss of generality that our algorithm~$\A$, on an input~$\sigma$, outputs an integer~$e$ which is a code for a partial computable function $\varphi_e$ from  $\N$ to $\mathcal{I}$ (our set of rational points in~$\Mcs$, described above) that is defined on the entire $\N$ or at some initial segment of $\N$, and $\rho(\varphi_e(n),\varphi_e(n+1)) < 2^{-n-1}$ when both $\varphi_e(n)$ and $\varphi_e(n+1)$ are defined. When this sequence is total, it converges to a measure $\mu$ with computable speed: $\rho(\varphi_e(n),\mu) \leq 2^{-n}$. 

This is not guaranteed by the definition of BD-success, but we may ``trim'' the algorithm by ensuring that indeed the sequence $\varphi_e(0), \varphi_e(1),\ldots$, whether finite or infinite, contains elements of $\mathcal{I}$ and satisfies the distance conditions where defined (by waiting until the conditions are checked; note that we have a strict inequality which will manifest itself at some moment, if true).

 Suppose now that for some index $e$, we do not know whether $\varphi_e$ is total, but we see that $\varphi_e(n)$ is defined for some $n$, and $\dfc(X \md \nu) > d$ holds for some $X$ and for \emph{all} measures $\nu$ at distance $\leq 2^{-n}$ of $\varphi_e(n)$. Then we already know, \emph{should $\varphi_e$ be total and converge to some $\mu$}, that $\dfc(X \md \mu) > d$. Thus we use the following notation: if $\A(\sigma)$ returns~$e$, then $\dfc(X \md \A(\sigma))$ is the quantity
\[
\sup \{d \md \exists n\, \varphi_e(n) \downarrow \text{ and } \dfc(X \md \nu) > d \text{ for all $\nu$ at distance $\leq 2^{-n}$ from $\varphi_e(n)$} \}.
\]
The supremum of an empty set (that appears, for example, if $\A(\sigma)$ is a code of an empty sequence) is considered to be $0$. Our function $\dfc(X \md \A(\sigma))$ has two key properties which are essential for the rest of our proof:
\begin{itemize}
\item[(a)] If $\A(\sigma)=e$, and $\varphi_e$ is total and converges to $\mu$, then $\dfc(X \md \A(\sigma)) = \dfc(X \md \mu)$. 
\item[(b)] $\dfc(X \md \A(\sigma))$ is lower semi-computable, uniformly in~$(X, \sigma)$.
\end{itemize}

Let us first prove that property (a) holds. Assuming $\A(\sigma)=e$, and $\varphi_e$ is total and converging to $\mu$, we know that $\mu$ is at distance $\leq 2^{-n}$ of $\varphi_e(n)$ for all~$n$. Thus, in the definition of $\dfc(X \md  \A(\sigma))$ every member of the set of $d$ on the right-hand side is at most $\dfc(X \md \mu)$. Thus $\dfc(X \md  \A(\sigma)) \leq \dfc(X \md \mu)$. On the other hand, if $\dfc(X \md \mu) > d$, by lower semicontinuity of the function $\dfc$, there is  $n$ such that $\rho(\nu,\mu) \leq 2^{-n}$ implies $\dfc(X \md \nu) > d$, and thus $\dfc(X \md  \A(\sigma)) > d$. Property (a) is proven. 

For property~(b), we use the fact that the space $\Mcs$ is effectively compact (one can effectively enumerate all covers of $\Mcs$ consisting of a finite union of open rational balls), together with the fact that the infimum of a lower semicomputable function on an effectively compact set is lower semicomputable, uniformly in a code for the effectively compact set (see~\cite{Gacs2005} for both facts). Thus, the predicate ``$\dfc(X \md \nu) > d \text{ for all $\nu$ at distance $\leq 2^{-n}$ from $\varphi_e(n)$}$" is computably enumerable uniformly in~$e,n,X,d$ (or said otherwise, the set of $(e,n,X,d)$ satisfying this property is an effectively open subset of $\N \times \N \times \cs \times \N$). Property~(b) follows. \\

Thanks to property~(a), when $\A(\sigma)=e$ and $\varphi_e$ is total and converges to a measure~$\mu$, we can safely identify $\A(\sigma)$ with $\mu$ and write $\A(\sigma)=\mu$. Additionally, we say that ``$\A(\sigma)$ is a measure" when $\A(\sigma)=\mu$ for some measure~$\mu$. \\

Now, for every pair of integers $(N,d)$, we define the set
\[
\WRONG(N,d) = \left\{ X \big| (\exists n \geq N) \; \dfc(X \md \A(X \uh n)) > d \right\}.
\]

This is the set of sequences $X$ on which the algorithm is ``visibly wrong" at some prefix of length $n \geq N$, for the deficiency level $d$. 

Note that $\WRONG(N,d)$, understood in this way, is effectively open uniformly in $(N,d)$ and is non-increasing in each of its parameters. The intersection of sets $\WRONG(N,d)$ for all $N$ and $d$ is some set $\WRONG$; as the name says, the algorithm $\A$ cannot BD-succeed on any sequence in this set. (Note that other reasons for failure are possible, e.g., $\A$ may not provide a measure on prefixes of some sequence.)

It is technically convenient to combine the two parameters $N$ and $d$ into one (even they are of different nature) and consider a decreasing sequence of sets $\WRONG(N)=\WRONG(N,N)$ whose intersection is $\WRONG$. 

We also consider a set $\SUCC(N,d)$ of all sequences $X$ such that $\A$ BD-succeeds on $X$ at level $N$ with deficiency $d$, i.e., 
\[
\SUCC(N,d) = \left\{ X \colon (\forall n \geq N) \, \left[ \text{$\A(X\uh n)$ is a measure, }  \dfc(X \md \A(X \uh n)) \le d \right]\right\}.
\]
The set $\SUCC(N,d)$ is a closed set. Indeed, it is an intersection of sets indexed by $n$, so we need to show that each of them is closed. For each $n$ there are finitely many possible prefixes of length $n$, so the first condition (``$\A(X\uh n)$ is a measure'') defines a clopen set. The second condition defines an effectively closed subset in each cylinder where $\A(X\uh n)$ is a measure. (Note that we do \emph{not} claim that $\SUCC(N,d)$ is \emph{effectively} closed, since the condition ``to be a measure'' is only a $\Pi_2$-condition.) By definition, the set $\SUCC(N,d)$ does not intersect the set $\WRONG(N,d)$.

The set $\SUCC(N,d)$ increases as $N$ or $d$ increase; the union of these sets is the set of all $X$ where $\A$ BD-succeeds; we denote it by $\SUCC$. Again we may combine the parameters and consider an increasing sequence of sets $\SUCC(N)=\SUCC(N,N)$ whose union in $\SUCC$.

All these considerations deal with the space of sequences. Now we switch to the space of measures and the class $\M$. We look what are the measures of sets $\WRONG(N)$ according to different $\mu\in\M$. Consider some threshold $x\in[0,1]$. There are two possible cases:
\begin{itemize}
\item there exist some number $N$, and some non-empty open set $\U \subseteq \M$ such that $\mu(\WRONG(N))\le x$ for all $\mu \in \U$. 
\item for every $N$ the set of points $\mu\in \M$ where $\mu(\WRONG(N))>x$ is dense in~$\M$.
\end{itemize}
There is some threshold where we switch from one case to the other, so let use take close values of $p<q$ (i.e., we take the difference $q-p$ to be much smaller than $\delta$ from the statement of the theorem; it would be enough to require that $q-p<\delta/10$) such that the first case happens for $q$ and the second one happens for $p$.

Choose some $N$ and an open ball $\B_0$ that has a non-empty intersection with $\M$ such that $\mu(\WRONG(N))\le q$ for all $\mu\in\B_0\cap\M$ (this is possible since the first case happens for $q$). 

\begin{lemma}
There exists a computable measure $\mu^*\in \B_0\cap \M$ such that\\ $\mu^*(\WRONG)\ge p$.
\end{lemma}

\begin{proof}
Since the second case happens for $p$, we can find some $\mu\in \B_0\cap \M$ such that $\mu(\WRONG(0))>p$. Since $\WRONG(0)$ is open, the same is true for some its clopen subset $C$, i.e., $\mu(C)>p$.  Note that $\mu(C)$ is a continuous function of $\mu$ for fixed clopen $C$, so we can find a smaller ball $\B_1\subseteq \B_0$ intersecting $\M$ such the $\mu(\WRONG(0))\ge \mu(C)>p$ for all $\mu\in \overline\B_1\cap\M$. Then, repeating the same argument, we find an even smaller ball $\B_2\subseteq \B_1$ intersecting $\M$ such that $\mu(\WRONG(1))>p$ for all $\mu\in \overline\B_2\cap\M$, then some $\B_3\subseteq \B_2$ such that $\mu(\WRONG(2))>p$ for all $\mu\in\overline\B_3\cap\M$, etc. Using the completeness of the space of measures, consider the intersection point $\mu^*$ of all $\B_i$ (we may assume that their radii converge to $0$ and that $\overline\B_{i+1}\subseteq\B_i$, and this guarantees the existence and the uniqueness of the intersection point). We have $\mu^*(\WRONG(i))>p$ for all $i$ (but $\mu^*(\WRONG(N))\le q$; the same in true for all subsequent sets $\WRONG(i)$ for all $i\ge N$). The continuity property for measure $\mu^*$ then guarantees that $\mu^*(\WRONG)\ge~p$.

Refining this argument, we can get a \emph{computable} measure $\mu^*$ with this property. Indeed, we may choose $\B_{i+1}$ in such a way that even the closed ball $\overline\B_{i+1}$ of the same radius is contained in $\B_i$; this property is enumerable. ``To have a non-empty intersection with $\M$'' is also an enumerable property (by assumption), and ``$\mu(\WRONG(i))>p$ for all $\mu\in\overline\B_{i+1}$'' is also an enumerable property (we may assume without loss of generality that $p$ is rational, and $\WRONG(i)$ is effectively open uniformly in $i$). So we can perform a search until $\B_{i+1}$ is found, and the sequence of $\B_i$ is computable, so $\mu^*$ is computable.\qed
\end{proof}

Now the argument goes as follows. Since $\mu^*$ is computable, the set $\SUCC$ should have $\mu^*$-probability at least $\delta$ by assumption. Success means that (at least) some measures are provided by the learning algorithm $\A$ for prefixes of sufficiently large length $M$. There are finitely many possible prefixes, and they correspond to finitely many computable measures $\mu_1,\ldots,\mu_s$. Then we choose a measure $\mu'$ orthogonal to all these measures and very close to $\mu^*$. We get the contradiction showing that  $\mu'(\WRONG(N))$ is almost $p+\delta$ (or more) and therefore exceeds $q$ which is not possible due to the choice of $\B_0$. To get the $\delta$-increase we use the fact that sequences that are $\mu'$-random cannot be $\mu_i$-random and should therefore have infinite deficiency. Let us now explain this argument in details.

Recall that we have chosen $N$ in such a way that $\mu(\WRONG(N))\le q$ for all $\mu$ sufficiently close to $\mu^*$. On the other hand, $\mu^*(\WRONG(M))\ge \mu^*(\WRONG)\ge p$ for all $M$.

Since $\mu^*(\SUCC)\ge\delta$, the continuity property of measures guarantees that $\mu^*(\SUCC(M))\gtrsim \delta$ for sufficiently large $M$, where $\gtrsim$ means inequality up to an additive error term that is very small compared to $\delta$ (in fact, $\delta/10$ would be small enough; we do not add more than $10$ inequalities of this type). Fix some $M$ that is large enough and greater than $N$ from the previous paragraph.

The set $\WRONG(M)$ is open and has $\mu^*$-measure at least $p$. Therefore, there exist a clopen set $C\subseteq \WRONG(M)$ such that $\mu^*(C)\gtrsim p$. Since the set $C$ is clopen, there exists some $K\ge M$ such that the $K$-bit prefix determines whether a sequence belongs to $C$ (the granularity of $C$ is at most $K$).

Now the Cantor space is split into $2^K$ intervals that correspond to different prefixes of length $K$. Some of these intervals form the set $C$ (and belong to $\WRONG(M)$ entirely). Among the rest, we distinguish good and bad intervals; good intervals correspond to prefixes for which the learning algorithm $\A$ produces a measure (whatever this measure is).  Let $\mu_1,\ldots,\mu_s$ be all measures that are produced by $\A$ for all good intervals (we have at most $2^K$ of them).

Note that $\SUCC(M)$ is covered by the good intervals. Indeed, it is disjoint with $\WRONG(M)$ and therefore with $C$, and also is disjoint with bad intervals by definition (since $K\ge M$, the algorithm $\A$ should produce a measure when applied to $K$-bit prefix). 

Now consider a measure $\mu'$ that is very close to $\mu$ and orthogonal to all $\mu_i$. (Our assumption allows us to get a measure very close to $\mu$ and orthogonal to a given computable measure; now we have several measure $\mu_1,\ldots,\mu_s$ instead of one, but this does not matter since we may consider their average: any measure orthogonal to the average is orthogonal to all $\mu_i$.)

Since the $\mu^*$-measure of $\SUCC(M)$ is almost $\delta$ (or more), and it is covered by good intervals, then $\mu^*$-measure of the union of good intervals is also almost $\delta$ (or more). The same is true for every measure $\mu'$ sufficiently close to $\mu^*$ since the union of good intervals is a clopen set.

No $\mu'$-random sequences can be $\mu_i$-random since the measures are orthogonal. This implies infinite deficiency, so all $\mu'$-random sequences in good intervals belong to $\WRONG(M)$. So the $\mu'$-measure of the part of $\WRONG(M)$ outside $C$ is almost $\delta$ (or more), and the part of $\WRONG(M)$ inside $C$ has $\mu'$-measure almost $p$ or more (this was true for $\mu$, and $\mu'$ is close to $\mu$). Together we get lower bound close to $p+\delta$ for $\mu'(\WRONG(M))$. And this gives us a contradiction, since $\mu'(\WRONG(M))\le \mu'(\WRONG(N))$, and the latter should be at most $q$ for all $\mu'$ close to $\mu$.  (Recall that we have chosen $q-p$ much smaller than $\delta$.)

This contradiction finishes the proof of Theorem~\ref{thm:main-strong}.

\section{Removing the deficiency boundedness requirement:\\ a positive result}\label{sec:positive}

We have established in Theorem~\ref{thm:main-weak} that, unsurprisingly, there is no total algorithm $\A$ that BD-succeeds on all sequences~$X$ that are random with respect to some computable probability measure. After proving this theorem, the authors conjectured that one could even prove the same result for UD-success. But it turns out that the situation is drastically different for UD-learning: we will show in this section that there \emph{is} a uniform learning algorithm in this model. 

\begin{theorem}\label{thm:positive-result}
There exists a total algorithm~$\A$ that UD-succeeds on every $X$ that is Martin-L\"of random with respect to some computable probability measure.
\end{theorem}

Recall that this means that for large enough $n$, $\A(X \uh n)$ is a (code of a) measure with respect to which~$X$ is random. However, ({\em a}) $\A(X \uh n)$ may be different for different values of~$n$, and ({\em b}) the randomness deficiency of $X$ with respect to $\A(X \uh n)$ is unbounded.

The proof of this result is inspired by a result of Harrington (reported in~\cite[Theorem 3.10]{CaseS1983} or \cite[Theorem VII.5.55, p.~139]{OdifreddiVolume2}) which states that there exists an algorithm to learn --- in the classical sense --- all computable sequences \emph{up to finitely many errors}. More precisely, there is a total algorithm~$\A$ such that for every computable sequence~$X$, for almost all~$n$, $\A(X \uh n)$ is a program for an almost everywhere defined function that differs from $X$ only in finitely many places. Indeed, let $\A(\sigma)$ be the program that, given input~$m$, spends time~$m$ searching for the minimal program computing some extension of~$\sigma$ and then runs this program, if found, on~$m$ (and returns, say,~$0$ if no such program is found). Let $e$ be the minimal program that computes~$X$. All smaller programs fail to compute~$X$ on some $k$ (by either being undefined or giving a wrong answer).  If $n$ is greater than all these $k$, then none of the smaller programs (than $e$) would qualify as a candidate for any $m$, and for large enough $m$ the program $e$ will be approved. (Note that $\A(\sigma)$ may be a non-total program even if $\sigma$ is long: we know only that it is defined on large enough values of~$m$.)

\begin{proof}
We will use an argument somewhat similar to Harrington's  to prove Theorem~\ref{thm:positive-result}. In this section, it is more convenient to consider measures as functions by identifying $\mu$ with the function $\sigma \mapsto \mu(\sigma)$ (here and in the rest of the section, $\mu(\sigma)$ is the abbreviation of $\mu([\sigma])$). 

It is also more convenient  to use an alternative characterization of Martin-L\"of randomness, via the Schnorr--Levin theorem, which states that if $\mu$ is a computable measure, a sequence~$X$ is Martin-L\"of random with respect to $\mu$ if and only if the prefix complexity of its prefixes is big:
\[
(\exists c)(\forall n) ~ \K(X \uh n) > -\log \mu(X \uh n) - c 
\]
(see for example~\cite{LiV2008}). We say that measure $\mu$ is \emph{exactly computable} when the function $\sigma \mapsto \mu(\sigma)$ is rational-valued and computable as a function from~$\fs$ to~$\Q$.  Of course, not all computable measures are exactly computable, but the following fact holds:

\begin{lemma}
If $X \in \cs$ is random with respect to a computable probability measure $\mu$, it is random with respect to some exactly computable probability measure~$\nu$. Moreover, one can suppose that $\nu(\sigma) >0$ for all strings~$\sigma$.
\end{lemma}

See~\cite{JuedesL1995} for a proof (essentially we approximate the given computable measure with enough precision using positive rational numbers). 

This lemma is convenient because it is possible to effectively list the family~$\F$ of partial computable functions~$\mu$ from $\fs$ to $\Q$  such that

\begin{itemize}
\item $\mu(\Lambda)=1$;
\item for every~$n$, either $\mu(\sigma)$ is defined for all strings $\sigma$ of length~$n$, or is undefined for all strings $\sigma$ of length~$n$;
\item if $\mu(\sigma 0)$ and $\mu(\sigma 1)$ are defined, $\mu(\sigma)$ is defined and is equal to $\mu(\sigma 0) +\mu(\sigma 1)$;
\item $\mu(\sigma)>0$ for all $\sigma$ on which $\mu$ is defined.
\end{itemize}

Let $(\mu_e)_{e \in \N}$ be an effective enumeration of all the functions in~$\F$. It is among these functions that, given a sequence~$X$, we are going to look for the `best candidate' $\mu_e$ such that $\mu_e$ is a measure (i.e., is total) and $X$ is random relative to~$\mu_e$. Suppose we are given a prefix $\sigma$ of~$X$. What is a good candidate $\mu_e$ for this $\sigma$? For this, we use the same approach as algorithmic statistics: a good explanation $\mu_e$ for a string $\sigma$ should (\emph{a}) be defined on $\sigma$, (\emph{b}) be simple, which is measured by the prefix complexity $\K(e)$, and (\emph{c}) give $\sigma$ a small `local' randomness deficiency, which we measure by the quantity $\ld(e,\sigma)=\max_{\tau \leqx \sigma} [-\log \mu_e(\tau) - \K(\tau)]$, with the convention that $\ld(e,\sigma)=\infty$ when $\mu_e(\tau)$ is undefined for some prefix $\tau$ of~$\sigma$. The Schnorr--Levin theorem mentioned above now can be reformulated as follows: the value
\[
d(X\md\mu_e)=\sup_{\tau\leqx X}[-\log \mu_e(\tau) - \K(\tau)]=\lim_n \ld(e,X\uh n)=\sup_n \ld(e,X\uh n)
\]
is finite if and only  if $\mu_e$ is a measure and $X$ is Martin-L\"of random with respect to $\mu_e$. In fact, $d(X\md e)$ is a version of randomness deficiency; for a fixed measure $\mu_e$ this quantity is equal to the deficiency~$\dfc$ of the previous section up to logarithmic precision (see, e.g.,~\cite{gacs-long} for details). 

Returning to algorithmic statistics, we combine the two quantities into a score function
           $$\score(e,\sigma) = \K(e) + \lceil \ld(e,\sigma) \rceil,$$
(as in golf, `score' is meant in a negative sense: a high $\score(e,\sigma)$ means that $\mu_e$ is not a good candidate for being a measure with respect to which~$\sigma$ looks random). Finally, we define the function $\BEST$ such that $\BEST(\sigma)$ is the value of $e$  that minimizes $\score(e,\sigma)$  (if there are several, we let $\BEST(\sigma)$ be the smallest one). That is,
$$\BEST(\sigma)=\min\{e\md (\forall e')\ \score(e,\sigma)\leq \score(e',\sigma) \}.$$

The first thing to observe is that $\BEST$ is computable in the halting set~$\mathbf{0}'$. Indeed, to compute $\BEST(\sigma)$, one can first find~$e$ such that $s=\score(e,\sigma) < \infty$ (this can be done computably). Then, using~$\mathbf{0}'$, one can find~$N$ such that $\K(e)> s$ (and thus $\score(e,\sigma)>s$) for all $e>N$. Finally, take $\BEST(\sigma)$ to be the number $e$ in $[0,N]$ that minimizes $\score(e,\sigma)$ (again taking the smallest one if there are several), which can be done effectively relative to~$\mathbf{0}'$ because $\score$ is itself computable relative to~$\mathbf{0}'$.

The core of the proof of Theorem~\ref{thm:positive-result} is the following lemma, which is of independent interest. It implies that learning measures in the EX sense, which we showed in the previous section to be impossible, becomes possible if one is given access to oracle $\mathbf{0}'$.

\begin{lemma}\label{lem:cvg-of-best}
Let $X$ be a sequence that is random with respect to some computable probability measure. The sequence of integers $\BEST(X \uh n)$ converges to a single value~$e^*$ such that $\mu_{e^*}$ is a measure, and $X$ is random with respect to $\mu_{e^*}$.
\end{lemma}

\begin{proof}
Fix such a sequence~$X$. For each~$e$, the sequence $\score(e,X \uh n)$ is nondecreasing and takes its values in $\N \cup \{\infty\}$, thus converges to some $S(e) \in \N \cup \{\infty\}$ . As we have said, the Schnorr--Levin theorem guarantess that $S(e) < \infty$ if and only if  $\mu_e$ is a measure and $X$ is Martin-L\"of random with respect to $\mu_e$. Thus we know that $S(e) < \infty$ for some~$e$ by our assumption that $X$ is Martin-L\"of random with respect to some computable probability measure. 

Let~$e^*$ be the index such that $S(e^*)$ is minimal among all~$S(e)$ (the smallest one if there are several). For any $i$ such that $K(i) > S(e^*)$, we have for any~$n$:
\[
\score(i,X \uh n) > S(e^*) \geq \score(e^*, X \uh n)
\] 
Thus $\BEST(X \uh n) \not= i$ for any~$n$. In other words, only the $j$ such that $K(j) \leq S(e^*)$  matter when selecting the best candidate for the sequence~$X$. Those~$j$ form a finite set. For all such~$j$, we know that $\score(j, X \uh n)$ is non-decreasing and eventually reaches its final value. After that, for all sufficiently large $n$, we have $\BEST(X\uh n)=e^*$.\qed
\end{proof}

At this point, we have seen that the function $\BEST$ does achieve the learning of measures we want, but unfortunately this function is only $\mathbf{0}'$-computable. By Schoenfield's limit lemma, this means that there exists a computable procedure which, given $\sigma$, generates a sequence $e_0, e_1, \ldots $ of integers that converges to $e_\infty = \BEST(\sigma)$. There is in general no way to compute $\mu_{e_\infty}$ from this sequence. However, what we \emph{can} do is combine all the $\mu_{e_i}$ together (being cautious about the fact that some $\mu_{e_i}$ may be partial) into a single computable measure $\nu$ such that $\nu > c \mu_{e_\infty}$ for some $c>0$, and this, by the Schnorr--Levin theorem, guarantees that everything that is random with respect to $\mu_{e_\infty}$, is also $\nu$-random.

More precisely, we have the following lemma.

\begin{lemma}\label{lem:gluing}
Let $f: \fs \rightarrow \N$ be a total $\mathbf{0}'$-computable function such that $\mu_{f(\sigma)}$ is a measure for all $\sigma$. Then there is a computable function $g$ such that $\mu_{g(\sigma)}$ is a measure for all~$\sigma$, and $\mu_{g(\sigma)} \geq c_\sigma \mu_{f(\sigma)}$ for some positive~$c_\sigma$. 
\end{lemma}

\begin{proof}
Consider the following effective procedure. On input $\sigma$, we use the Schoenfield limit lemma to effectively get a sequence $e_i$ converging to $e_\infty = f(\sigma)$. Initially all $e_i$ are considered ``candidates". We then apply a filtering process that deletes some of these candidates. Recall that the corresponding $\mu_{e_i}$ are elements of $\mathcal{F}$. We compute in parallel all $\mu_{e_i}(\tau)$ for all pairs $(i,\tau)$ for which $e_i$ is still a candidate. If we find two candidates $e_i, e_j$ and $\tau$ such that $\mu_{e_i}(\tau)$ and $\mu_{e_j}(\tau)$ are both defined and not equal different from each other, then we remove $e_i$ and $e_j$ from the list of candidates. This way we ensure, since the sequence converges, that from some point on, for any candidate $e_i$,  the corresponding function $\mu_{e_i}$ is equal to $\mu_{e_\infty}$ on its domain (but $\mu_{e_i}$ may be partial). Indeed, each bad candidate (i.e., an $e_i$ such that $e_i \not= e_\infty$) may destroy at most one good candidate, and by assumption almost all candidates are good.

Now we let $\mu_{g(\sigma)}$ to be a computable measure $\nu$ constructed in the following way. First, let $\nu(\Lambda)=1$. Then we compute the conditional probabilities $\nu(x0)/\nu(x)$ and $\nu(x1)/\nu(x)$ level by level. When computing them on level $N$, we use for the computation the conditional probabilities for some candidate that remains alive after $N$ steps of filtering process. (Any of them could be used, for example, we may take the one with smallest computation time. Note the at least one good candidate remains, so we will not wait forever.)

As we have seen, starting from some level only good candidates remain, so the conditional probabilities above this level are the same for $\mu_{f(\sigma)}$ and $\nu$. Since by assumption all values of all measures are positive, this guarantees the required inequality. \qed 
\end{proof}

We can now put all pieces together to prove Theorem~\ref{thm:positive-result}. Applying the previous lemma to $f=\BEST$, we have a computable function $g$ such that for every~$\sigma$, the measure $\mu_{g(\sigma)}$ dominates, up to a multiplicative constant, the measure $\mu_{\BEST(\sigma)}$. For every~$X$ that is random with respect to some computable measure, we know, by Lemma~\ref{lem:cvg-of-best}, that $\mu_{\BEST(X \uh n)}$ is eventually constant and equal to a measure with respect to which $X$ is random. This measure is dominated (up to multiplicative constant) by $\mu_{g(X \uh n)}$, thus $X$ is also random with respect to $\mu_{g(X \uh n)}$ (change in the measure increases the deficiency at most by $O(1)$). This finishes the proof. \qed
\end{proof}

\noindent \textbf{Acknowledgements}. Laurent Bienvenu and Santiago Figueira acknowledge the support of the Laboratoire International Associ\'e ``INFINIS''. Laurent Bienvenu and Alexander Shen also acknowledge the support of ANR-15-CE40-0016-01 RaCAF grant.




%
\bibliographystyle{alpha} 
\bibliography{bfms-learning}

\end{document}